\documentclass[11pt]{article}
\usepackage{epsfig}
\usepackage[utf8]{inputenc} 

\usepackage{rotating}
\usepackage{graphicx}
\usepackage{verbatim}

\usepackage{url}
\usepackage{color}
\usepackage{pstricks}
\usepackage{fancyvrb}
\usepackage{makeidx}

\VerbatimFootnotes

\begin{document}

\newcommand{\bm}[1]{\mbox{\boldmath$#1$}}

\def\mvec#1{{\bm{#1}}}   % vector in math mode; bold

\title{Probability, propensity 
and \\ probabilities of propensities 
(and of probabilities)\footnote{Invited 
contribution to the proceedings MaxEnt 2016 
based on the talk
given at the workshop (Ghent, Belgium, 10-15 July 2016), 
supplemented by work done 
within the program  Probability and Statistics in Forensic Science
at the Isaac Newton Institute for Mathematical Sciences, Cambridge,
under the EPSRC grant EP/K032208/1.
}
}

\author{G.~D'Agostini \\
Universit\`a ``La Sapienza'' and INFN, Roma, Italia \\
{\small (giulio.dagostini@roma1.infn.it,
 \url{http://www.roma1.infn.it/~dagos})}
}

\date{}

\maketitle

%\vspace{-0.5cm}
\begin{abstract}
The process of doing Science in condition of
uncertainty is illustrated with a toy experiment
in which the inferential and the forecasting aspects
are both present. The fundamental aspects of probabilistic
reasoning, also relevant in real life applications, 
 arise quite naturally and the resulting discussion
among non-ideologized, free-minded people offers an opportunity
for clarifications.
\end{abstract}

%\vspace{-0.3cm}
{\small 
\begin{flushright}
{\sl \small ``I am a Bayesian in data analysis,}\\
{\sl \small   a frequentist in Physics''}\\
(A PhD student in Rome, 2011)\\
\mbox{} \\
{\sl  \small ``You see, a question has arisen,}\\
{\sl \small about which we cannot come to an agreement,}\\
{\sl \small probably because we have read too many books''}\\
(Brecht's Galileo) \\ 
\mbox{} \\
{\sl \small  ``The theory of probabilities is basically}\\ 
{\sl \small               just common sense reduced to calculus''}\\
(Laplace)%
\end{flushright}
}
\mbox{}
\vspace{-18cm}
{\small 
\begin{flushleft}
\tt %DESY 95-242 \hfill ISSN 0418-9833 \newline %( questo numero e' di DESY )
Isaac Newton Institute \\
for Mathematical Sciences\newline
Cambridge, UK \newline 
INI preprint NI16052 (December 2016)
\end{flushleft}
}
\vspace{15.5cm}

\vspace{-1.3cm}
\section{Introduction}
Much has been said and written about probability. Therefore,
instead of presenting the different views, or accounting for
its historical developments, I go straight
to an example, which I like to present
as an experiment, as indeed it is: the boxes and the
balls are real and they represent the `Physical World'
about which we `do Science,' that is 1) we try,
{\em somehow}, to gain 
our knowledge about it by making observations; 2) we try, {\em somehow},
to anticipate the results of future observations.
`Somehow' because we usually start
and often remain in conditions of uncertainty.
So, instead of starting by saying
``probability is defined as such and such'',
I introduce the toy experiment, explain the rules of the `game,'
clarifying what can be directly observed and what can only be guessed,
and then let   the discussion go, guiding it with proper questions
and helping it by evaluating interactively numbers of interest (some lines of
R code are reported in the paper for the benefit of the reader).
Later, I make the `players' aware of the implications of their answers
and choices.
And even
though initially some of the numbers do not come out right
 -- the example is simple enough  that
rational people will finally agree on the numbers
of interest -- the main concepts do:
subjective probability as degree of belief;
physical `probability' as propensity of systems to behave in a given way;
the fact that we can be uncertain about the values of propensity,
and then assign them probabilities; and even that degrees of beliefs
can  themselves be
uncertain and often expressed in fuzzy terms like `low', 'high',
`very high' and so on --  when this is the case they
need to be defuzzified before they can
be properly used within probability theory, without the need to invent
something fancy in order to handle them.
Other % important 
points touched in the paper are the myth that
propensities are only related to long-term relative frequencies
and the question of verifiability of events
subject to probabilistic assessments. 

\section{Which box? Which ball?}
\begin{figure}
\centerline{\includegraphics[width=0.7\linewidth]{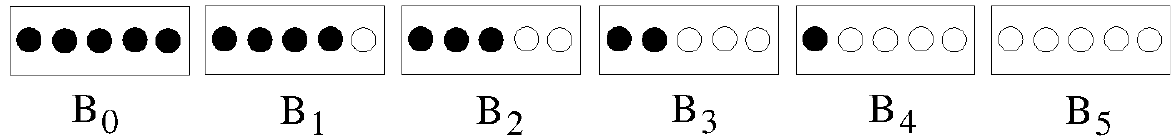}}
\mbox{}\vspace{-0.7cm}
\caption{\small \sf A sketch of the six boxes of the toy experiment.
The index refers to the number of white balls.}
\label{fig:sixbox}
\end{figure}
The `game' begins by showing 
six boxes (Fig.~\ref{fig:sixbox}), 
each containing five balls.\footnote{Those who understand
Italian might form an idea of a real session watching
a video of a conference for the general public organized by
the University of Roma 3 in June 2016
(\url{http://orientamento.matfis.uniroma3.it/fisincittastorico.php#dagostini})
and available on YouTube
(\url{https://www.youtube.com/watch?v=YrsP-h2uVU4}).
}
One box has only black balls, another four Black and one White,
and so on.
One box, hereafter `B$_?$', is taken at random out of the six
and we start the game.
At each stage, we have to guess which box has been chosen
and what color ball will be selected in a random extraction.
We then extract a ball, observe its color and replace it 
into the box \cite{AJP}.

From the point of view of measurements,
the uncertain number of white balls plays the role of the
value of a physical quantity; the two colors
the possible empirical observations.
The fact that we deal with
a discrete and small set of possibilities,
both for the `measurand' and the empirical `data',
only helps in clarifying the reasoning.
Moreover, one of the rules of the game is that
we are forbidden to look
inside the box, in the same way that we cannot open an electron and
read its mass and charge in a hypothetical label inside it.

\subsection{Initial situation}
At the first stage the answers to the questions are prompt
and unanimous:
we consider all boxes equally likely,
thus assigning 1/6 probability to each of them;
we consider Black and White
      equally likely too, with probabilities 1/2.
Not satisfied with these answers, I also encourage `players' to express
their confidence on the hypotheses of interest
by means of a virtual lottery at zero entry cost.\footnote{For this
purpose this kind
of lotteries are preferable to normal bets, although hypothetical
and even those with small amount
of money (value and amount of money are well known for not being
proportional), in order to allow people to freely choose what they
consider more credible, without incurring
the so called {\em loss aversion bias}.}
Specifically, I ask, 
if you are promised a large prize
for making the correct prediction,
which box or ball color would you
choose?
More precisely, is there any reason
to prefer a particular color or a particular composition?
Also, in this case there is a general consensus on the fact that
any choice is equally good, in the sense that there is no reason
to be blamed if we finally miss the prize.

\subsection{An intriguing dilemma: B$_?$ Vs  B$_{E}$}
At this point a new box $B_E$ with equal number of black
and white balls
is shown to the audience. In contrast with B$_?$, everyone 
can now check its content
(the box I actually use  contains 5 White and 5 Black).
In this case we are only uncertain about
the result of picking a ball,
and, again, everyone considers Black and White equally probable.

Then a new virtual lottery is proposed, with
a prize associated to the extraction of  {\em White} from
either box. Is it preferable to
choose B$_?$ or B$_{E}$? That is, is there any special reason to
opt for either box?
This time
the answer is not always unanimous and depends on the audience.
Scientists, including PhD students,
{\em tend to} consider -- but there are practically always exceptions! --
the outcomes equally probable and therefore they say there
is no rational reason to prefer either box.
But in other contexts, including
seminars to people who have jobs of high responsibility,
there is a sizable proportion, often the majority,
of those who definitely prefer B$_{E}$ (and, by the way,
they had already stated, or accepted without objections, that
Black and White were equally likely also from this box!)

The fun starts in case (practically always) when there are people
in the audience having shown
a strong preference in favor of B$_{E}$, and later I change
the winning color. For example, I say, just for the sake of entertainment, that the prize in case of White was supposed to be
offered by the host of the seminar. But since I prefer black,
as I am usually dressed that way,
\underline{I} will pay for the prize, but attaching it to {\em Black}.
As you can guess,
those who showed  indifference between B$_?$ and  B$_{E}$
keep their opinion (and stare at me in a puzzled way).
But, curiously, also those who had previously chosen with
full conviction  B$_{E}$
stick to it.
The behavior of the latter is quite irrational
(I can understand one can have strange reasons to consider White more likely from $B_E$, but for
the same reason he/she should consider Black more likely from
$B_?$) but so common that it even has a name,  the {\em Ellsberg paradox}.
(Fortunately the kind of people attending my seminars
repent quite soon, because they are easily
convinced -- this is the simplest explanation --
that, after all,
the initial situation with $B_?$ is absolutely equivalent to
an extraction at random out of 15 Black and 15 White, the fact that
the 30 balls are clustered in boxes being irrelevant.)

\subsection{Changing our mind in the light of the observations}
Putting aside box B$_{E}$,
from which there is little to learn for the moment,
we proceed with our `measurements'
on box B$_?$.
Imagine now that the first extraction gives {\em White}.
There is little doubt that the
observation {\em has to}\footnote{In this particular case it is
clear that `it has to', but in general `it might'. See for example
footnote \ref{fn:status_information} and pay attention
that conditional probabilities
might be not intuitive and a formal guidance is advised.
}
change somehow
our confidence on the box composition
and on the color that will result from the next extraction.

As far as the box composition is concerned,
 $B_0$ is ruled out, since ``this box cannot give
white balls,'' or, as I suggest,
 ``this cause cannot produce the observed effect.''
In other words, hypothesis $B_0$ is `falsified,'
i.e. the probability {\em we assign to it} 
drops instantly to zero.
But what happens to the others? The answer of the large
majority of people, with remarkable exceptions
(typically senior scientists),
is that the other compositions remain equally likely,
with probability values then rising from 1/6 to 1/5.

The qualitative answer to the second question is basically correct,
in the sense that it goes into the right direction:
the extraction of
 White {\em becomes} more probable,\footnote{Please compare this expression,
``the extraction of White {\em becomes} more probable'', with 
``the probability {\em we assign to it}'', used above. The former should 
be, more correctly, ``we assign higher probability 
to the extraction of White'', as it will be clear later. 
For sake of conciseness and avoiding pedantry, in this paper 
I will often use imprecise expressions of this kind,
as used in every day language.} 
``because $B_0$ has been ruled out.''
But, unfortunately, the quantitative answer never comes out right,
at least initially.
In fact, at most, people say that the probability of White
rises to 15/25, that is 3/5, or 60\%, just from the arithmetic
of the remaining balls after $B_0$ has been removed from
the space of possibilities.

The answers ``remaining compositions equally likely''
and ``3/5 probability of  White'' are both {\em wrong}, but
they are at least consistent, the second being a
logical consequence of the first, as can easily be shown.
Therefore, we only need to understand what is wrong with the first answer,
and this can be done at a qualitative level, just with
a bit of hand waving.\footnote{See
e.g. \url{https://www.youtube.com/watch?v=YrsP-h2uVU4}
from 48:00
(in Italian). }
Imagine the hypothetical case of a long sequence of
White, for example 20, 50 or even 100 times
(I remind that extractions are followed by re-introduction).
After many observations we start to be highly confident that
we are dealing with box $B_5$, and therefore the probability of
White in a subsequent extraction
approaches unity. In other words, we would be {\em highly
surprised} to extract a black ball, already after 20 White in a row,
not to speak after 50 or 100, although we do not consider such an event
absolutely impossible. It is simply highly improbable.

It is self-evident that, if after many observations we reach such a situation
of {\em practical certainty}, then every extraction has to contribute
a little bit. Or, differently stated, each observation
has to provide a bit of evidence
in favor of the compositions with larger proportions of white balls.
And, therefore, even the very first observation
has to break our symmetric state of uncertainty over
the possible compositions.
How? At this point of the discussion there is a kind of general
enlightenment in the audience:
the probability has to be proportional to the number
of white balls of each hypothetical composition, because
``{\em boxes with a larger proportion of white balls tend to produce more
easily White},'' and therefore ``{\em White comes easier from $B_5$ than
$B_4$, and so on}.''

\section{Updating rules}

\subsection{Updating rule for the ``probabilities of the causes''} 
The heuristic rule resulting from the discussion is
\begin{eqnarray}
P(B_?=B_i\,|\,\mbox{W},I) & \propto & \pi_i\,,
\end{eqnarray}
where $ \pi_i=i/N$, with $N$ the total number of balls in box $i$,
is the white ball proportion and
$I$ stands for
all other available information regarding the experiment.
[In the sequel we shall use the shorter notation
$P(B_i\,|\,\mbox{W},I)$ in place of $P(B_?=B_i\,|\,\mbox{W},I)$,
keeping instead always explicit the `background' condition $I$.]
But, since  the probability
$P(\mbox{W}\,|\,B_i,I)$ of getting White from box $B_i$
is trivially  $\pi_i$ (we shall come back to the reason) we get
\begin{eqnarray}
P(B_i\,|\,\mbox{W},I) & \propto & P(\mbox{W}\,|\,B_i,I)\,.
\end{eqnarray}
This rule is obviously not general, but depends on the fact
that we initially considered all boxes equally likely,
or $P(B_i\,|\,I) \propto 1$, a convenient notation in place
of the customary $P(B_i\,|\,I) = k$, since common factors are irrelevant.
So a reasonable {\em ansatz} for the updating rule,
consistent with the result of the discussion, is
\begin{eqnarray}
P(B_i\,|\,\mbox{W},I) & \propto &   P(\mbox{W}\,|\,B_i,I) \cdot
P(B_i\,|\,I)\,.
\end{eqnarray}
But if this is the proper updating rule, it has to hold after
the second extraction too, i.e. when $P(B_i\,|\,I)$ is replaced by
$P(B_i\,|\,\mbox{W},I)$, which we rewrite as
$P(B_i\,|\,\mbox{W}^{(1)},I)$ to make it clear that such a
probability depends {\em also} on the
observation of White in the first extraction. We have then
 \begin{eqnarray}
P(B_i\,|\,\mbox{W}^{(1)},\mbox{W}^{(2)},I) & \propto &
P(\mbox{W}^{(2)}\,|\,B_i)\cdot
                 P(B_i\,|\,\mbox{W}^{(1)},I) \,,
\end{eqnarray}
and so on. By symmetry, the updating rule in case Black (`B')
were observed is
 \begin{eqnarray}
P(B_i\,|\,\mbox{B},I) & \propto & P(\mbox{B}\,|\,B_i)\cdot P(B_i\,|\,I)\,,
\end{eqnarray}
with $P(\mbox{B}\,|\,B_i) = 1-\pi_i$. After a sequence of $n$ White
we get therefore $P(B_i\,|\,\mbox{`$n$W'},I) \propto \pi_i\,^n$. For example
after 20 White we are -- we must be! --
98.9\% confident to have chosen $B_5$ and 1.1\%
$B_4$, with the remaining possibilities `practically'
ruled out.\footnote{Here is the result with a single line of R code:\\
\verb|> N=5; n=20; i=0:N; pii=i/N; pii^n/sum(pii^n)| \\
\verb|[1] 0.000000e+00 1.036587e-14 1.086940e-08 3.614356e-05 1.139740e-02 9.885665e-01| \\
(And, by the way, this is a good example of the importance
of a formal guidance in assessing probabilities: according to my experience,
after a sequence of 5-6 White,
people are misguided by intuition and tend to believe
box $B_5$ much more than they rationally \underline{should}.)
}

If we observe, continuing the extractions,
 a sequence of $x$ White
and $(n-x)$ Black we get\footnote{
Here is the R code for the example of 20 extractions resulting in 5 White:\\
\verb|> N=5; n=20; i=0:N; pii=i/N; x=5; pii^x * (1-pii)^(n-x) / sum( pii^x * (1-pii)^(n-x) )| \\
\verb|[1] 0.000000e+00 6.968411e-01 2.979907e-01 5.167614e-03 6.645594e-07 0.000000e+00| \\
(Note how using this code we can focus on the essence of what it is going
on, instead of being `distracted' by the math of  the normalization.)
}
 \begin{eqnarray}
P(B_i\,|\,n,x,I) & \propto & \pi_i^x\,\left(1-\pi_i\right)^{n-x}\,.
\label{eq:P(Bi_i|seq)}
\end{eqnarray}
But, since there is  a one-to-one relation between $B_i$ and $\pi_i$,
we can write
\begin{eqnarray}
P(\pi_i\,|\,n,x,I) & \propto & \pi_i^x\,\left(1-\pi_i\right)^{n-x},
\label{eq:P(pi_i)}
\end{eqnarray}
an apparently `innocent' expression on which we shall comment later.

\subsection{Laplace's `Bayes rule'}
As a matter of fact, the above updating rule can be shown to result
from probability theory, and I find it magnificently
described in simple words by Laplace in what he calls
``{\sl the fundamental principle
of that branch of the analysis of chance
that consists of reasoning {\it a posteriori} from events
to causes}''\,\cite{Laplace}:\footnote{In the light of
Brecht's quote by Galileo you might be surprised to find quite
some quotes in this paper. But there are books and books.
}
{\small
\begin{quote}
      {\sl
      ``The greater the probability of an observed event given any one
of a number of causes to which that event may be attributed,
the greater the likelihood
% \,\footnote{Note that
%here likelihood is the same as probability, and has nothing
%to do with what statisticians call `likelihood' -- reading
%directly the original French version might help,
%also taking into account
%that two hundred years ago the nouns were not as specialized
%as they now are.}
 of that cause \{given that event\}.
The probability of the existence of any one of these causes
 \{given the event\} is thus a fraction
whose numerator is the probability of the event given the cause,
and whose denominator is the sum of similar probabilities,
summed over all causes. If the various causes are not equally probable
{\em a priori}, it is necessary, instead of the probability of the event
given each cause, to use the product of this probability
and the possibility
of the cause itself.''\,\cite{Laplace}
      }
\end{quote}
}
\noindent
Thus, indicating by $E$ the effect and by $C_i$ the $i$-th cause,
and neglecting normalization,  Laplace's
{\em fundamental principle} is as simple as
\begin{eqnarray}
P(C_i\,|\,E,I) &\propto& P(E\,|\,C_i,I) \cdot P(C_i\,|\,I)\,,
\label{eq:Bayes1}
\end{eqnarray}
from which we learn a simple rule that teaches us how
to update the ratio of probabilities we assign to
two generic causes $C_i$
and $C_j$ (not necessarily mutually exclusive):
\begin{eqnarray}
\frac{P(C_i\,|\,E,I)}{P(C_j\,|\,E,I)} &=&
\frac{P(E\,|\,C_i,I)}{P(E\,|\,C_j,I)}
\cdot \frac{P(C_i\,|\,I)}{P(C_j\,|\,I)}\,.
\end{eqnarray}
Equation (\ref{eq:Bayes1}) is a convenient way to
express the so-called {\em Bayes rule} (or `theorem'), while the
last one shows explicitly how the ratio of the probabilities of two causes
is updated by the piece of evidence $E$ via the so called
{\em Bayes factor} (or {\em Bayes-Turing factor}\,\cite{IJ_Good}).
Note the important implication 
of Equation (\ref{eq:Bayes1}): we cannot update the probability of a cause, 
unless it becomes {\em strictly} falsified, if we not consider
at least another fully specified cause \cite{BadMath,WavesSigmas}. 

\subsection{Updating the probability of the next observation}
Coming to the probability of White in the second extraction,
it is now clear why $15/25=3/5=60\%$ is wrong: it assumed
the remaining five boxes equally
likely,\footnote{~This would have been the correct answer to a different
question: probability of White from a box taken at random
among boxes $B_{1-5}$,
that is $B_?^{(1-5)}$. Ruling out $B_0$ by hand at the very
beginning is quite different from ruling it out as a consequence
of the described experiment. The status of information is different
in the two cases and also the resulting probabilities will usually be
different! [Please note that a different state of information {\em might}
change probability, but not necessarily it does.
For example $P(\mbox{W}^{(1)}\,|\,I) =
P(\mbox{W}^{(11)}\,|\,5\mbox{B},5\mbox{W},I)$
just by symmetry.
Conditioning is subtle!]
\label{fn:status_information}
}
while they are not.
Also in this case maieutics helps: it becomes suddenly
clear that we have to assign a higher `weight' to the compositions
we consider more likely. That is, in general and remembering that
the weights $P(B_i\,|\,I)$ sum up to unity,
\begin{eqnarray}
P(\mbox{W}\,|\,I) &=& \sum_i P(\mbox{W}\,|\,B_i,I)\cdot
                                 P(B_i\,|\,I)\,.
\label{eq:P(W|all_B)}
\end{eqnarray}
After the observation of White in the first
extraction we then get\footnote{~Here is the numerical result obtained with R:\\
\verb|> N=5; i=0:N; pii=i/N; ( PBi = pii/sum(pii) ); sum( pii * PBi ) | \\
\verb|[1] 0.00000000 0.06666667 0.13333333 0.20000000 0.26666667 0.33333333| \\
\verb|[1] 0.7333333|
}
\begin{eqnarray}
P(\mbox{W}^{(2)}\,|\,\mbox{W}^{(1)},I) &=&
\sum_i P(\mbox{W}^{(2)}\,|\,B_i,\mbox{W}^{(1)},I)\cdot
                                 P(B_i\,|\,\mbox{W}^{(1)},I) \nonumber \\
                 &=&
\sum_i P(\mbox{W}\,|\,B_i,I)\cdot
                                 P(B_i\,|\,\mbox{W}^{(1)},I)\,,
\label{eq:P(W2|W1)}
\end{eqnarray}
where $P(\mbox{W}^{(2)}\,|\,B_i,\mbox{W}^{(1)},I)$ has been rewritten
as $P(\mbox{W}\,|\,B_i,I)$ since, assuming a particular composition,
the probability of White is the same in every extraction.
Moreover,
since $\pi_i = P(\mbox{W}\,|\,B_i)$,
we can rewrite Equation (\ref{eq:P(W2|W1)}), in analogy with
Equation (\ref{eq:P(pi_i)}), i.e. replacing $B_i$ by $\pi_i$,
as \begin{eqnarray}
P(\mbox{W}^{(2)}\,|\,\mbox{W}^{(1)},I) &=&
\sum_i \pi_i\cdot P(\pi_i\,|\,\mbox{W}^{(1)},I)\,,
\label{eq:PW2_pi}
\end{eqnarray}
which will deserve comments later.

\section{Where is probability?}
The most important outcome of the discussion related
to the toy experiment is in my opinion that,
although people do not immediately get the correct numbers,
they find it quite natural that relevant changes
of the available information
have to modify somehow the probability of the box composition
and of the color resulting in a future extraction,
although {\em the box remains the same},
i.e. nothing changes inside it.\footnote{Curiously,
for strict frequentists
the probability that $B_?$ contains $i$ white balls
makes no sense because, they say, either it does or it doesn't.
}
Therefore the crucial,
rhetorical  question follows: {\em Where is the probability?}
Certainly \underline{not in the box}!

At this point, as a corollary, it follows that,
if someone just enters the room and does not
know the result of the extraction, he/she will reply to
our initial questions
 exactly as we initially did.  In other words,
there is no doubt that the probability has to depend
on the {\it subject} who evaluates it, or
{\small
 \begin{quote}
      {\sl
      ``Since the knowledge may be different with different persons
      or with the same person at different times, they may anticipate
      the same event with more or less confidence, and thus different numerical
      probabilities may be attached to the same event}.''\,\cite{Schrodinger}
\end{quote}
}
\noindent
      If follows that probability is always conditional
      probability, in the sense that
{\small
\begin{quote}
      {\sl
      ``Thus whenever we speak loosely of `the probability of an event,'
      it is always to be understood: probability with regard to a certain
      given state of knowledge.''\,\cite{Schrodinger}
      }
\end{quote}
}
\noindent
     So, more precisely, $p=P(E)$ should always be understood as $p=P(E\,|\,I_S(t))$,
     where $I_S(t)$ stands for the information available
     to the subject $S$ who evaluates $p$
     at time $t$.\footnote{~The notation used above is consistent
     with this statement, in the sense that the conditions
     appearing in $P(B_i\,|\,I)$, $P(B_i\,|\,\mbox{W}^{(1)},I)$ and
     $P(B_i\,|\,\mbox{W}^{(1)},\mbox{W}^{(2)},I)$
     can be seen seen as $I_S(t)$ evolving with time.}
     It is disappointing that many confuse `subjective' with `arbitrary',
     and they are usually the same who make use of arbitrary
     formulae not based on
     probability theory, that is the {\em logic of uncertainty},
     but because they are supported by the Authority Principle,
     pretending they are `objective'.\footnote{It is curious to remark
     that there
     are, or at least there were,
     also Bayesians `afraid' of subjective probability\,\cite{GdA_ME2000}.}

\section{What is probability?}
  A third quote by Schr\"odinger summarizes the first two
and clarifies what we are talking about:
{\small
\begin{quote}
      ``{\sl
      Given the state of our knowledge about everything that could
possibly have any bearing on the coming true. . . the numerical
probability $p$ of this event is to be a real number by the
indication of which we try in some cases to setup a
quantitative measure of the strength of our conjecture or
anticipation, founded on the said knowledge, that the event
comes true
      \,\cite{Schrodinger}
      }
\end{quote}
}
\noindent
Probability {\em is not} just ``a number between 0 and 1 that
satisfies some basic rules'' (`the axioms'),
as we sometimes hear and read,
because such a `definition' says
nothing about what we are talking about.
If we can understand probability statements
it is because we are able, so to say,
to map them in some `categories' of our mind,
as we do with space and time (although for
values far from those we can feel  directly
with our senses
we need some means of comparison, as when we say
``30 times the mass of the sun'', and rely on numbers).

Think for example of two generic events $E_1$ and $E_2$ such that
$p_1=P(E_1\,|\,I)$ and
$p_2=P(E_2\,|\,I)$. Imagine also that we
have our reasons -- either we have evaluated the numbers, or
we trust  somebody's else evaluations --
to believe that $p_1$ is  {\em much} larger that $p_2$,\footnote{~Note
also this very last statement, to which we shall return
at the end of the paper.}
 where  `much' is added in order to make
our {\em feeling} stronger.  It is then a matter of fact that:
 ``the strength of our conjecture'' strongly favors $E_1$;
 we expect (``anticipate'')
      $E_1$ much more than $E_2$;
 we will be highly surprised if  $E_2$ occurs, instead
of $E_1$.\footnote{~As a real example, in my talk
at MaxEnt 2016 I analyzed the football match France-Portugal,
played right on the first day of the workshop, so that
everybody (interested in football) had fresh in their minds
the reaction of fans of the two teams, as shown on TV,
and also that of people in pubs in Ghent
(slides are available at
\url{http://www.roma1.infn.it/~dagos/prob+stat.html#MaxEnt16_2}).
}
Or, in simpler words, {\em we believe $E_1$ to occur much more than $E_2$}.

\subsection{Ideas, beliefs and probability}     
In other terms, finally calling things with their name,
we are talking about {\em degree of belief},
and references to the deep and thorough analysis
of David Hume are deserved.
The reason we can communicate with
each other our degrees of belief (``I believe this more than that'')
is that our mind understands what we are talking
about, although\footnote{What Hume says about
probability reminds me of the famous reflection by Augustine of Hippo
about time: {\sl ``Quid est ergo tempus? Si nemo ex me quaerat, scio; si quaerenti explicare velim, nescio.`` -- ``What then is time? If no one asks me, I know what it is. If I wish to explain it to him who asks, I do not know.''}
(\url{https://en.wikiquote.org/wiki/Augustine_of_Hippo}.)
Indeed, as a creature living in a
hypothetical Flatland
has no intuition of
how a 3D world would be, so a hypothetical intelligent humanoid
 `determinoid,' living in a (very boring) world
in which all phenomena happen with extreme regularity, would
have not developed the concept of probability.
}
{\small
\begin{quote}
      {\sl
      ``This operation of the mind,
which forms the belief of any matter of
fact, seems hitherto to have been one of the greatest mysteries of
philosophy\\
\ldots \\
When I would explain \{it\}, I
scarce find any word that fully answers the case, but am obliged to have
recourse to every one's feeling, in order to give him a perfect notion
of this operation of the mind.''
\,\cite{Hume_Treatise}%
}.
\end{quote}
}
\noindent
In fact, since {\sl ``nothing is more free than the
imagination of man''}\,\cite{Hume_Enquiry}, we can conceive
all sorts of ideas,
just combining others. But we do not consider them all
believable, or equally believable:
 {\sl ``An idea assented to {\it feels} different
from a fictitious idea, that the fancy alone presents to us: And this
different feeling I endeavour to explain by calling it a
superior {\it force},
or {\it vivacity}, or {\it solidity},
or {\it firmness}, or steadiness.''}\,\cite{Hume_Treatise}
(italics original.)

An easy evaluation is when we have a set of {\em equiprobable}
cases, a proportion of which leads to the event of interest
(neglect for a moment the first sentence of the quote):
{\small
\begin{quote}
      {\sl
       [``Though there be no such thing as {\em Chance} in the world; our
ignorance of the real cause of any event has the same influence on the
understanding, and begets a like species of belief or opinion.''] \\
      ``There is certainly a probability, which arises from a superiority of
chances on any side; and according as this superiority encreases, and
surpasses the opposite chances, the probability receives a
proportionable encrease, and begets still a higher degree of belief or
assent to that side, in which we discover the superiority. If a dye were
marked with one figure or number of spots on four sides, and with
another figure or number of spots on the two remaining sides, it would
be more probable, that the former would turn up than the
latter.''\,\cite{Hume_Enquiry}%
}.
\end{quote}
}
\noindent
This is the reasoning we use to assert that the probability of White
from box $B_i$ is proportional to $i$, viz. $P(\mbox{W}\,|\,B_i,I)=\pi_i$.
Instead, the precise reasoning which allows us to evaluate
the probability of White from $B_?$ in the light of the previous extraction
was not discussed by Hume (for that we have to wait until
Bayes \cite{Bayes},
and Laplace for a thorough analysis \cite{Laplace_ATP}),
but the concept of probability still
holds. For example, after four consecutive white balls the probability of
White in a fifth extraction becomes about 90\%.
That is, assuming the calculation has been done correctly,
we are essentially so confident to extract White from $B_?$
as we would from a box containing 9 white
balls and 1 black.\footnote{~The exact number of
$P(\mbox{W}^{(5)}\,|\,4\mbox{W},I)$ is 90.4\%, as it can be easily checked
with R:\\
\verb|> N=5; n=4; i=0:N; pii=i/N; ( PBi=pii^n/sum(pii^n) ); sum(pii * PBi)| \\
\verb|[1] 0.00000000 0.00102145 0.01634321 0.08273749 0.26149132 0.63840654| \\
\verb|[1] 0.9039837|
}

\section{Physical probability? }
Going back to the previous quote by Hume, an interesting, long debated
issue is whether there is  {\sl ``such a thing as Chance in the world''},
or if, instead, probability arises {\em only} because
of {\sl ``our ignorance of the real cause of any
event.''}\footnote{The second
position, popularized by Einstein's ``God does not play dice'',
is related to the so-called Laplace Demon,
``An intellect which at a certain moment would know all forces
that set nature in motion, and all positions of all items of
which nature is composed, if this intellect were also vast
enough to submit these data to analysis,
it would embrace in a single formula the movements
of the greatest bodies of the universe and those of the tiniest atom;
for such an intellect nothing would be uncertain and the future just
like the past would be present before its eyes.''~\cite{Laplace}
}
This is a great question which I like to tackle in a very
pragmatic way, re-wording the first sentence of the quote:
whatever your opinion might be,
{\sl ``the influence on the understanding''} is the same.
If you assign 64\% probability to event $E_1$ and 21\% probability to
$E_2$ (and 15\% that something else will occur)
you simply believe (and hence your mind ``anticipates'')
$E_1$ much more that $E_2$, no matter what  $E_1$ and of $E_2$
refer to,
provided you are {\em confident on the probability values}
(please take note of this last expression).

For example, the events could be White and Black from a box
containing 100 balls, 64 of which White, 21 Black,
and the remaining of other colors.
But $E_1$ could as well be the decay of the `sub-nuclear' particle
K$^+$ into a {\em muon} and a {\em neutrino},  and
$E_2$ the decay of the same particle
into  two {\em pions}
(one charged and one neutral).\footnote{The {\em branching ratios}
of K$^+$ into the two `channels' are
$\mbox{BR}(\mbox{K}^+ \rightarrow  \mu^+\nu_\mu) = (63.56\pm 0.11)\%$
and
$\mbox{BR}(\mbox{K}^+ \rightarrow  \pi^+\pi^0) =
(20.67\pm 0.08)\%$\,\cite{PDG}. \\
By the way, I do not think that Quantum Mechanics needs special rules
of probability. There the mysteries are related to the
weird properties of the wave function $\psi(x,t)$. Once you apply the
rules -- ``shut up and calculate!'' has been for long time
 the pragmatic imperative --
and get `probabilities'
(in this case `propensities', as we shall see) all the rest
is the same as when you calculate `physical probabilities' in other systems.
Take for example the brain-teasing
single photon double slit experiment
(see e.g. \url{https://www.youtube.com/watch?v=GzbKb59my3U}).
From a purely probabilistic
point of view the situation is quite simple. Applying the rules of
Quantum Mechanics,  if we open only slit $A$ we get
the pdf $f_A(x\,|\,A,I)$; if we open only $B$ we get
 $f_B(x\,|\,B,I)$; if we open both slits  we get
$f_{A\& B}(x\,|\,{A\& B},I)$. Why should
$f_{A\& B}(x\,|\,{A\& B},I)$ be just a superposition of
$f_A(x\,|\,A,I)$ and  $f_B(x\,|\,B,I)$? In fact within probability
theory there is no rule which relates them. We
need a  model to evaluate each of them
and the best we have are the rules of Quantum Mechanics. 
%A different story is to get an intuition of these rules. 
Once we have got the above pdf's all the rest follows
as with other common pdf's. In particular, if 
we get e.g. that  $f_A(x_1\,|\,A,I) >> f_{A\& B}(x_1\,|\,{A\& B},I)$
we believe that
a photon will be detected `around' $x_1$, if we open only slit $A$,  
much more than if we open both slits. And, similarly, if we plan to
repeat the experiment a large number of times, 
we expect to detect `many more' photons `around' $x_1$ if 
only slit $A$ is open than if both are.  
That's all.
A different story is to get an intuition of the rules of Quantum Mechanics.
}
Thus,
as  we consider the 64\% probability of the
K$^+$ to produce a muon and a neutrino a physical property
of the particle, similarly it can be {\em convenient}
to consider the  64\%
probability
of the box to produce white balls a physical property
of that box, in addition to its mass and dimensions.
(It is interesting to pay attention
to the long chain of somebody else's
beliefs, implicit when e.g. a physicist uses a published branching ratio to
form his/her own belief on the decay of a
particle.\footnote{I like, as historian Peter Galison
puts it: {\sl ``Experiments begin and end in a matrix of beliefs.
\ldots Beliefs in instrument type, in programs of experiment
enquiry, in the trained, individual judgments about every local behavior
of pieces of apparatus.''}\,\cite{Galison} Then beliefs are propagated
within the scientific community and then outside.
But,  as recognized,  methods
from `standard statistics' (first at all the infamous p-values)
tend to confuse even experts and spread unfounded beliefs
through the scientific community as well as
among the general public \cite{BadMath,WavesSigmas}, that in the meanwhile
is developing `antibodies'
and is beginning to mistrust striking scientific results and, I am afraid,
sooner or later also scientists and Science in general.
}
And something similar occurs for other quantities and in other domains
of science and in any other human activity.)

\subsection{Propensity vs probability}
Back to our toy experiment,
I then see no problem saying that box $B_i$
\underline{has} probability $\pi_i$ to produce white balls, meaning that
such a `probability' is a physical property of the box, something that
measures its {\em propensity} (or {\em bent}, {\em tendency},
{\em preference})\footnote{I have no strong preference on the name, and my
propensity in favor of `propensity' is because it is less used in
ordinary language (and despite the fact that this noun
is usually associated to Karl Popper,
an author I consider quite over-evaluated).}
to produce white balls.

It is a matter of fact that, if we have full confidence that
a {\em physical}\,\footnote{Note the extended meaning of `physical',
not strictly related to Physics, but to `matters of fact' of all kinds,
including for example biological, sociological or economic systems
{\em believed} to have propensities to behave in different ways.
} system
has propensity $\pi$ to produce event $E$,
then we shall use $\pi$ to form the
{\sl ``strength of our conjecture or anticipation''} of its occurrence,
that is $P(E\,|\,\pi,I)=\pi$.\footnote{I had heard
that this apparent obvious statement goes under the name
of Lewis' {\em Principal Principle}
(see e.g. \url{http://plato.stanford.edu/entries/probability-interpret/}).
Only at the late stage of writing this paper
I bothered to investigate a little
more about that `curious principle' and found out
Lewis' {\em Subjectivist's Guide
to Objective Chance}\,\cite{Lewis}, in which his very basic concepts,
outlined in a couple of dozen of lines at the beginning
of the article, are amazingly in tune
with several of the positions I maintain here.
\label{fn:Lewis}
}
 But it is often the case in real life that, even
if we hypothesize that such a propensity does exist,
we are not certain about its value, as it happens
with box $B_?$. In this case we have to take into account
all possible values of propensity. This is the meaning of
Equation (\ref{eq:PW2_pi}), which we can rewrite in more general terms as
\begin{eqnarray}
P(E\,|\,I) &=& \sum_i \pi_i\cdot P(\pi_i\,|\,I)\,.
\end{eqnarray}
We can extend the reasoning to a continuous
set of $\pi$, indicated by $p$ for its clear meaning of
the parameter of a Bernoulli process,
% (e.g. that appearing in the binomial distribution)
to which we associate then a
 probability density function (pdf), indicated  by
$f(p\,|\,I)$:\footnote{It becomes now clear the meaning 
of Equation (\ref{eq:P(pi_i)}), which we can rewrite
as
\begin{eqnarray*}
f(p\,|\,n,x,I) & \propto & p^x\,\left(1-p\right)^{n-x},
\end{eqnarray*}
having assumed a continuity of propensity values, and
having started our inference from a uniform {\em prior}, that is
$f(p\,|\,I) = 1$. \\
The normalized version of the above equation is
%\mbox{}\vspace{-0.3cm}
\begin{eqnarray*}
f(p\,|\,n,x,I) &  = & \frac{(n+1)!}{x!\,(n-x)!}\,p^x\,\left(1-p\right)^{n-x}.
\end{eqnarray*}
\mbox{}\vspace{-0.4cm}
\label{fn:f(p)_norm}
}
\begin{eqnarray}
P(E\,|\,I) &=&  \int_0^1\! p\,f(p\,|\,I)\,\mbox{d} p
\label{eq:PW2_p}
\end{eqnarray}
The special case in which our {\em probability}, meant as {\em degree of
belief}, coincides with a particular value of {\em propensity},
is when $P(\pi_i\,|\,I)$ is 1 for a particular $i$, or $f(p\,|\,I)$ is a Dirac delta-function.
This is the difference between boxes $B_?$ and $B_E$.
In $B_E$ our degree of belief of $1/2$ on White or Black is directly
related to its assumed propensity 
to give balls of either colors. In $B_?$ a numerically identical degree of
belief arises from averaging all possible propensity values
(initially equally likely). And therefore the
{\sl ``strength of our conjecture or anticipation''}~\cite{Schrodinger}
is the same in the two cases.
Instead, if we had at the very beginning only the boxes
with at least one white ball,
the probability of White from $B_?^{(1-5)}$
becomes, applying the above formula,
$\sum_{i=1}^{5}(i/5)\times (1/5) = 3/5$.
%(If, instead, $B_0$ was ruled out as a result of the extraction
%of a white ball, then the different compositions, and related
%propensities, are differently likely, and we cannot use any longer
%the common $1/5$ weight, as we have seen above.)

We are clearly talking about {\em probabilities of
propensities}, as when we are interested in
detector {\em efficiencies}, or in {\em branching ratios}
of unstable particles (or in the proportion of the population
in a country that shares a given character or opinion, or the
many other cases in which we use a binomial distribution,
whose parameter $p$ has, or might be given, the meaning of propensity).
But there are other cases in which probability has no propensity
interpretation, as in the case of the probability of a box composition,
or, more generally, when we make {\em inference} on the
{\em parameter of a model}. This occurs for
instance in our toy experiment %, sticking to our toy experiment,
when we were talking about $P(B_i\,|\,I)$, a concept to which
no serious scientist objects, as well as he/she has nothing against
talking e.g. of $90\%$ probability that the
mass of a black hole lies within a given interval of values
(with the exception of a minority of
highly ideologized guys).

\subsection{Probability, propensity and (relative)  frequency}
A curious myth is that physical probability, or propensity,
has ``only a frequentist interpretation'' (and therefore ``physicists
must be frequentist'', as
ingenuously stated by the Roman  PhD student quoted on the first page).
But it seems to me to be more a question of education, based on the
dominant school of statistics in the past century
(and presently),\footnote{Here
is, for example, what David Lewis (see Footnote
\ref{fn:Lewis}) writes in Ref.~\cite{Lewis} (italics original):
{\sl ``Carnap did well to distinguish two concepts of probability,
insisting that both were legitimate and useful
and that neither was at fault because it was not the other.
I do not think Carnap chose quite the right two concepts,
however. In place of his `degree of confirmation', I would put
{\it credence} or {\it degree of belief}; in place
of his `relative frequency in the long run', I would put {\it chance}
or {\it propension}, understood as making sense in the single case.''}
More or less what I concluded when I tried to read Carnap
about twenty years ago: his first choice means nothing (or at least
it has little to do with probability); the second
does not hold, as I am arguing here.
}
rather than a real logical necessity.

It is a matter of fact that ({\em relative})
frequency and probability
are somehow connected within probability theory, without the need
for identifying the two concepts.
%\footnote{{\bf Expression probabilities
%with natural frequencies: confusione fra proporzioni e frequenze,
%e frazioni\ldots } }
\begin{itemize}
\item A {\em future} frequency $f_n$ in $n$ {\em independent}
      `situations' (not necessarily `trials'),\footnote{To make it clear,
      what is important to is that $p$ is ({\em about}) the same, 
      and that our assessments are independent. It does not matter
      if, instead, the events have a different meaning, like e.g. 
      tails tossing a coin, odd number rolling a die, and so on. 
      The emphasized `about' is because $p$ itself could be 
      uncertain, as we shall see later. In this case we need to 
      evaluate the expectation of $f_n$ taking into account the 
      uncertainty about $p$. 
      }
      to each of which we assign probability $p$, has expected value
      $p$ and `standard uncertainty' decreasing with increasing $n$ as
      $1/\!\sqrt{n}$,
      though all values 0, $1/n$, $2/n$, \ldots, 1 are
      {\em possible} (!).
      This a simple result of probability theory, directly related to the
      binomial distribution, that goes under the name of Bernoulli's theorem,
      often misunderstood with a  `limit', in
      the calculus's sense. Indeed $f_n$ {\em does not} ``tend to'' $p$,
      but it is simply {\em highly improbable}
      to observe $f_n$ far from $p$,
      for large values of $n$.\footnote{Related to this there is the
      usual confusion between a probability distribution and a distribution
      of frequencies. Take for example
      a quantity that can come in many possibilities, like
      in a binomial distribution with $n=10$ and $p=1/2$.
      We can think of repeating the trials a large number
      of times and then, applying Bernoulli's theorem
      to each of the eleven possibilities, we consider it very unlikely
      to observe values of relative frequencies in each `bin'
      different from the probabilities evaluated from the binomial
      distribution. This is why we highly expect -- and we
      shall be highly surprised at the contrary! -- a frequency distribution
      (`histograms') very similar in shape to the
      probability distribution, as you can easily `check' playing with \\
      \verb|n=10000; x=rbinom(n, 10, 0.5); barplot(table(x)/n, col='cyan')|\\
      \verb|barplot(dbinom(0:10,10,0.5), col=rgb(1,0,0,alpha=0.3), add=TRUE)| \\
       That's all! Nothing to do
      with the ``frequency interpretation of probability'',
      or with the ``empirical law of Chance''
      (see Footnote \ref{fn:empirical_laws}).}
      In particular,
      under the assumption that a system
      has a constant propensity $p$ in a large number of trials,
      we shall consider very ``unlikely to observe $f_n$
      far from $p$.''\footnote{Obviously, if you
      make an experiment of this kind, tossing
      regular coins or dice a large number of times,
      you will easily find relative frequencies of a given face
      around 1/2 or 1/6, respectively 
       as simulated with this 
      line of R: \\
      \verb|p=1/2;  n=10^5; sum( rbinom(n, 1, p) ) / n| \\
      But it is just because, in the Gaussian large number approximation,
      $P(|f_n-1/2| > 1/\sqrt{n}) = 4.6\%$,
       and therefore
      $f_n$ will
      {\em usually} occur around $1/2$ [although all $(n+1)$ values 
      between 0 and 1 are possible, 
      with probabilities $P(f_n=x/n) = 2^{-n}n! \left(x!(n-x)!\right)^{-1}$]. 
      Not because
      there is a kind of `law of nature'
      -- ``legge empirica del caso'', in Italian books,
      i.e.  ``empirical law of Chance'' --  `commanding' that
      {\em frequency has to tend
      to probability}, thus supporting the popular
      lore of late numbers at lotto hurrying up in order to obey it.
      In the scientific literature and in text books, not to speak
      about popularization books and article, it should be strictly
      forbidden to call `laws' the results of asymptotic theorems,
      because they can be easily misunderstood.
      \mbox{}[For example we read (visited 11/11/2016) in
      \url{https://it.wikipedia.org/wiki/Legge_dei_grandi_numeri}
      that ``the law of large numbers, also called empirical law of chance or Bernoulli's theorem [\ldots] describes \ldots''
      (total confusion! --
      see also \url{https://en.wikipedia.org/wiki/Law_of_large_numbers}
      and
      \url{https://en.wikipedia.org/wiki/Empirical_statistical_laws}).] 

      Moreover, it should be avoided to teach that e.g. probability 1/3
      means that something will occur to 1/3 of the elements of a 
      `reference class', {\em i}) first because a false sense of regularity 
      can be easily induced in simple minds, which will then complain
      that the ``the probabilities were wrong'' if no event of that
      kind occurred in 9 times; {\em ii}) second 
      because such `reference classes' might not exist, 
      and people should be trained in understanding degrees of belief
      referred to individual events. 
      \label{fn:empirical_laws}
      }  Reversing the reasoning, if we observe
      a given $f_n$ in a large number of trials, common sense
      suggests that the `true $p$' should lie not too far
      from it, and therefore our degree of belief
      in the occurrence of  a future event of that kind
      should be about $f_n$.
\item More precisely, the probability of a future event
      can be
      % formally
      mathematically
      related, under suitable assumptions,
      to the frequency of {\em analogous}\footnote{$E^{(1)}$
      is the success in the first trial, $E^{(2)}$
      the success in the second trial, and so on. Speaking
      about ``the realization of the same event'' is quite incorrect,
      because events $E^{(i)}$ are different. They can be at most
      analogous. We indicate here, instead, by $E$ the generic future
      event of the kind of $E^{(1)}$-$E^{(n)}$, i.e. for example
      $E=E^{(n+1)}$.}  events $E^{(i)}$ that occurred
      in the past.\footnote{~It is a matter of fact that,
because of evolution or whatever mechanism you might think about,
the human mind always looks for regularities.
This is how
 Hume puts it  (italics original):
 {\sl ``Where different effects have been found to
follow from causes, which are to {\it appearance} exactly similar, all these
various effects must occur to the mind in transferring the past to the
future, and enter into our consideration, when we determine the
probability of the event. Though we give the preference to that which
has been found most usual, and believe that this effect will exist, we
must not overlook the other effects, but must assign to each of them a
particular weight and authority, in proportion as we have found it to be
more or less frequent.''}\,\cite{Hume_Enquiry}
}     For example, assuming that a system has
      %unknown
      propensity $p$,
      after $x$ occurrences  (`successes')
      in $n$ trials
      we assign different beliefs to the different values of $p$
      according to a probability density function $f(p\,|x,n,I)$,
      whose
      % normalized
      expression has been reported in Footnote
      \ref{fn:f(p)_norm}.
      In order to take into account all possible values of $p$ we have
      to use Equation (\ref{eq:PW2_p}), in whose r.h.s. we recognize
      the {\em expected value} of $p$.
%      {\em the most probable
%      value has no  meaning, at least
%      talking about model parameters}).
      We get then the  famous Laplace
      {\em rule of succession} (and its limit for large $n$ and $x$),
      \begin{eqnarray}
      P(E\,|\,x,n,I) \!\!& = &\!\! 
\mbox{E}[p\,|\,x,n,I] =
\int_0^1\! p\,\frac{(n+1)!}{x!\,(n-x)!}\,p^x\left(1-p\right)^{n-x} \mbox{d}p = \frac{x+1}{n+2}
      \rightarrow \ \frac{x}{n} = f_n\,,
      \nonumber \\
      &&
      \label{eq:Laplace_rule_succession}
      \end{eqnarray}
      which can be interpreted as follows. If
      we {\em i}) consider the propensity of the system constant;
       {\em ii})  consider all values of $p$ a priori equally likely
                 (or the weaker condition of
                 all values between 0 and 1  possible,
                 if $n$ is `extraordinary large');
       {\em iii})  perform a `large' number of independent trials,
%      \begin{itemize}
%      \item[{\em i})] we consider the propensity of the system
%      constant;
%      \item[{\em ii})] we considered all values of $p$ a priori equally likely
%                 (or at least all value between 0 and 1 possible);
%      \item[{\em iii})]  we have performed a `large' number of independent trials,
%      \end{itemize}
       then
       the degree of belief we should assign to a future event
       is
       basically
       % essentially
       the observed past frequency.
       Equation (\ref{eq:Laplace_rule_succession})
       can then be seen as a mathematical proof
       that what the human mind does by intuition and
       ``custom'' (in Hume's sense) is  quite reasonable.
       But the formal guidance of probability theory makes clear the
      assumptions, as well as  the limitations of the result.
      For example, going back to our six box example, if after
      $n$ extractions we obtained $x$ White, one could be tempted
      to evaluate the probability of the next White from the observed
      frequency $f_n=x/n$,
      instead of, as probability theory teaches,
      firstly evaluating the probabilities
      of the various compositions from Equation (\ref{eq:P(Bi_i|seq)})
      and then the probability of White from (\ref{eq:P(W|all_B)}).
      The results will not be the same and the latter
      is amazingly `better'\footnote{
      To get an idea, repeat several times the following lines
      of R code which simulate \verb|n| extractions with re-introduction
      from box \verb|ri|,
      calculate the number of White,
      infer the probability of the box compositions,
      and finally evaluate the probability of a next White
      and compare it with the relative frequency.
      There is no miracle in the result, it is
      just that {\em the probabilistic formulae
      are  using all available information in the best possible way}: \\
      \verb|N=5; i=0:N; pii=i/N; ri=1; n=100; s=rbinom(n,1,pii[ri+1]); ( x=sum(s) )| \\
\verb|( PBi =  pii^x * (1-pii)^(n-x) / sum( pii^x * (1-pii)^(n-x) ) )| \\
\verb#cat(sprintf("P(W|sequence) = %.10f;  x/n = %.4f \n", sum( pii * PBi ), x/n))#
\label{fn:miracle}
      }\cite{AJP}.
\end{itemize}
There is another argument against the myth that
physical probability is `defined' via the long-term
frequency behavior. If propensity $p$ can be seen as a
parameter of a physical system, like a mass or the radius
of the sphere associated with the  shape of an object, then, as other
parameters, it might change with time too,
i.e. in general we deal with $p(t)$.
 It is then self-evident
that different observations will refer to propensities
at different times, and there is no way to get a long-term frequency at
a given time. At most we can make sparse measurements at different
times, which could still be useful, if we have a model
of how the propensity might change with
time.\footnote{I would like to make a related comment
on another myth concerning the scientific method, according to which
 ``replication is the cornerstone of Science''.
This implies that, if we take this principle literally, much of what
we nowadays consider
Science is in reality non-scientific
(can we repeat measurements concerning a particular supernova,
or two particular black holes merging with emission of
gravitational waves?).
And if you ask, they will tell you that this principle goes back
to none other than Galileo, who instead wrote\cite{Galileo}
% in his
% {\em Dialogue Concerning the Two Chief World Systems}
that
{\sl ``The knowledge of a single effect acquired
  by its causes opens our mind to understand
  and ensure us of other effects
  without the need of doing experiments''}
({\sl ``La cognizione d'un solo effetto acquistata per
            le sue cause ci apre l'intelletto a 'ntendere ed assicurarci
            d'altri effetti senza bisogno di ricorrere alle esperienze''}).
Doing Science is not just collecting
(large amounts of) data, but properly framing
them in a causal model of Knowledge.
}

\section{Abrupt end of the game -- Do we need verifiability?}
\label{sec:verifiability}
There is another interesting lesson that we can learn from our
six box toy experiment.\footnote{What is nice in this practical
session, instead of abstract speculations, is that the people
participating in the discussion have developed their degrees of beliefs,
and therefore, when the box is taken away, they cannot say that
what they were thinking (and feeling!) is not valid anymore.}
After some time the game has to come to an end, and the audience expects that I finally show the composition of box $B_?$.
Instead, I take it, put it
back together with the others and shuffle all them well. As you might imagine,
the reaction to this unexpected end is surprise and disappointment.
Disappointment because it is human to seek the `truth'.
Surprise because they didn't pay attention, or
perhaps didn't take me seriously,
when I said at the very beginning
that ``we are forbidden to look
inside the box, as we cannot open an electron and
read its mass and charge in a hypothetical label.'' %inside it.''

The reason for this unexpected ending of the game is twofold. First,
because scientists (especially students) have to learn, or to remember,
that when we make measurements we remain in most cases in a condition
of uncertainty.\footnote{See e.g. Feynman's quote 
at the end of the paper.} 
And not only in physics. Think, for example,
of forensics. How many times judges and jurors will finally know
with Certainty if the defendant was really guilty or
innocent?\footnote{If you worry about
 these issues, then you might be interested in the
Innocence Project, \url{http://www.innocenceproject.org/}.}
(We know by experience that we have to distrust even so-called
confessed criminals!)

The second reason is related to the question of the
{\em verifiability}
of the events about which we make probabilistic assessments.
Imagine, that during our toy experiment we made 6 extractions, getting
White twice, as for example in the following simulation in R.
(Note that if you run the lines of code as they are, deleting \verb|ri|
immediately after it is used in the second line,
you will never know the true composition! If you want to get
exactly the probability numbers of the last two outputs shown below,
without having to wait to get
\verb|x| equal  2, as it resulted here,
then just force its value.)\\
{\small
\verb|   > N=5; i=0:N; pii=i/N; n=6| \\
\verb|   > ri = sample(i, 1)|\\
\verb|   > ( s=rbinom(n,1,pii[ri+1]) ); rm(ri)| \\
\verb|   [1] 0 0 1 1 0 0| \\
\verb|   > ( x=sum(s) )   # nr of White|  \\
\verb|   [1] 2|  \\
\verb|   ( PBi =  pii^x * (1-pii)^(n-x) / sum( pii^x * (1-pii)^(n-x) ) ) | \\
\verb|   [1] 0.00000000 0.34594595 0.43783784 0.19459459 0.02162162 0.00000000| \\
\verb|   > sum( pii * PBi )| \\
\verb|   [1] 0.3783784| \\
}
At this point we have 44\% belief to have picked $B_2$
and only 2.2\% $B_4$; and 38\% belief to get White in a further
extraction. And these degrees of belief should be maintained, 
even if, afterwards,
we lose track of the box.\footnote{Note that many statements 
concerning scientific and historical `facts' are of this kind.}
This is like when we say that
 a plane {\em was} at a given instant in a given
cube of airspace with a given probability.
Or, more practically \cite{cep}, imagine you
are in a boat on the sea or on a lake, not too far from the shore,
so that you are able, e.g. using Whatsapp
on your smartphone,  to send to a friend
your GPS position, including its accuracy. The location is a point,
whose accuracy is defined by a radius such that
``there is a 68\% probability that
the true location is inside the circle.''\footnote{\url{https://developer.android.com/reference/android/location/Location.html#getAccuracy()}}
This is a statement that normal people, including experienced scientists,
understand and accept without problems
and which our mind uses to form a consequent degree of belief,
the same as when thinking of the probability of
a white ball being extracted blindly from a box
that contains 68 white and 32 black balls.
And practically nobody  has
concerns about the fact that {\em such an event cannot be verified}.
Exceptions are, to my knowledge, strict frequentists and strict definettians
(but I strongly doubt that they do not form in their mind
a similar degree of belief, although they cannot `professionally' admit it.)
In fact, for different reasons,
it is forbidden to scholars and practitioners of both schools
to talk about probability of hypotheses in the most general case, including
probability that true values are in a given
interval. For example neither of them could talk
of the probability that
the mass of Saturn is within a given interval,
as instead it was done by Laplace, to whom was perfectly
clear the hypothetical character of the so called
{\em coherent bet}.\footnote{Here is how
Laplace reported his uncertainty on value
of the mass of Saturn got by Alexis Bouvart:
%\begin{quote}
{\sl ``His [Bouvard] calculations give him
the mass of Saturn as 3,512th part of that
of the sun. Applying my
probabilistic formulae
 to these observations,
I find that the odds are 11,000 to 1 that the error in this
result is not a hundredth of its value.''} \cite{Laplace}
%\end{quote}
That is $P(3477 \le  M_{Sun}/M_{Sat}
\le 3547\,|\,I(\mbox{Laplace})) =  99.99\%\,.$
Note how the expression {\sl ``the odds \underline{are},''} indicates
he was talking of a fair bet, viz. a coherent bet.
Moreover it is self evident that
such a bet cannot be, strictly speaking, settled, but it rather had
an {\em hypothetical}, {\em normative} meaning.
(And Laplace was also well aware of the non linearity between
quantity of money and its `moral' value, so that
a bet with such high odds could never be agreed in practice
and it was just a strong way to state a probability.)
} As they would not accept talking about the most probable orbit
(``orbitam maxime probabilitatem''), or 
the probability
that a planet is at given point in the sky, as instead  
did Gauss when he derived his way the normal distribution
from the conditions (among others) that {\em i})  all points
were {\em a priori} equally likely ({\sl ``ante illas observationes [\ldots]
aeque probabilia fuisse''});  {\em ii}) the maximum 
of the {\em posterior} ({\sl ``post illas observationes''}) 
had to be equal to the arithmetic average of the observations
\cite{Gauss}.

\subsection{Probability of probabilities 
(and of odds and of Bayes factors)}
The issue of `probability of probability' has already been
discussed above, but in the particular case in which
the second `probability' of the expression was indeed
a propensity [and I would like to insist on the fact
that whoever is interested in probability
distributions of the Bernoulli parameter
$p$, that is in something of the kind $f(p\,|\,I)$,
is referring, explicitly or implicitly, to probabilities of propensities].
I would like now to move to the more general case,
i.e. when we refer to uncertainty about our degree of belief.
And, again, I like to approach the question in a pragmatic way,
beginning with some considerations.

The first is that we are often in situations in which we are reluctant
to assign a precise value to our degree of belief, because
``we don't know'' (this expression is commonly heard).
But if you ask ``is it then 10\%?'', the answer can be
``oh, not that low!'', or ``not so high!'' depending
on the event of interest. In fact it rarely occurs that we know absolutely nothing about the fact,\footnote{
{\sl ``If we were not ignorant there would be no probability,
 there could only be certainty. But our ignorance cannot
 be absolute, for then there would be no longer any probability
 at all. Thus the problems of probability may be classed
 according to the greater or less depth of our ignorance.''}
\cite{Poincare}
}
and in such a
case we are not even interested in evaluating probabilities
(why should we assign probabilities if we don't even know what
we are talking about?).
%In fact
%\begin{quote}
%{\sl ``If we were not ignorant there would be no probability,
% there could only be certainty. But our ignorance cannot
% be absolute, for then there would be no longer any probability
% at all. Thus the problems of probability may be classed
% according to the greater or less depth of our ignorance.''}
%\cite{Poincare}
%\end{quote}

The second is that the probability of probability,
in the most general sense, is already included
in the following, very familiar formula of probability theory,
valid if $H_i$ are all the elements of a complete class of hypotheses,
\begin{eqnarray}
P(A\,|\,I) &=& \sum_iP(A\,|\,H_i,I)\cdot P(H_i\,|\,I)\,.
\label{eq:weighted_average}
\end{eqnarray}
We only need the courage
to read it with an open mind:
Equation (\ref{eq:weighted_average}) is simply an average
of conditional probabilities, with weights equal to
probabilities of each contribution. But in order
to read it this way at least $P(H_i\,|\,I)$ must have the meaning
of degree of belief, while $P(A\,|\,H_i,I)$ can represent propensities or also
degrees of belief.

Probability of probabilities could refer
to evaluations of somebody's else 
probabilities,\footnote{Italians might be pleased to remember
Dante's ``Cred'io ch'ei credette ch'io credesse che \ldots'' (Inf. XIII, 25), 
expressing beliefs of beliefs of beliefs 
(``I believe he believed that I believed that\ldots''), 
roughly rendered in verses as ``He, as it seem'd, believ'd,
that I had thought [that]\ldots'' 
(\url{https://www.gutenberg.org/files/8789/8789-h/8789-h.htm#link13}).} 
as e.g. in game theory,
 but they are also important in all those
important cases of real life
in which direct assessments are done by experts
or when {\em sensitivity analysis} leads to a spectrum of possibilities.
For example, one might
evaluate his/her degree of belief around 80\%, but it could be as well,
perhaps with some reluctance,
75\% or 85\%, or even `pushed' down to 70\% or up to 90\%. With
suitable questions\footnote{For example we can ask the range
of virtual coherent bets one could accept
in either direction, or `calibrate' probabilistic judgements
against boxes with balls of different colors (or other mechanical
or graphical tools).}
it is possible then to have an idea of the range of possibilities,
in most cases with the different values not equally likely
(sharp edges are never reasonable).
For example, in this case it
could be a triangular distribution peaked at 80\%.
This way of modeling the uncertainties on degrees of belief
is similar to that recommended by the ISO's GUM
({\em Guide to the expression of uncertainty in measurement}\,\cite{ISO})
to model uncertainties due to systematic effects.
After we have modeled uncertain probabilities we can use the
formal rules of the theory to `integrate over' the possibilities,
analytically or by Monte Carlo
(and after some experience you might find out that, if you have several uncertain contributions, the details of the models
are not really crucial, as long as mean
and variance of the distributions
are `reasonable').
The only important remark is to be careful with probabilities
approaching 0 or 1. This can be done using log scale for {\em intensities of belief}, for the details of which I refer to
 \cite{Columbo}
 [in particular
Sections 2.4, 3.1, 3.3 and 3.4 (especially Footnote 22),
and Appendix E] and references therein.

Once we have broken the taboo of {\em freely} speaking
(because in reality we already somehow do it)
of probabilities of probabilities,
it is obvious that there is no problem to extend
this treatment of uncertainty to related quantities,
like odds and Bayes factor,
%of {\em probability of odds} and of
%{\em probability Bayes factors},
{\em i}) as a simple propagation from
uncertain probabilities;
{\em ii}) in direct assessments by experts. For example,
direct assessments of odds are currently performed for many
real-life events. Direct (`subjective')
assessments of Bayes factors were indeed
envisaged in Ref.~\cite{Columbo}.

\section{Conclusions}
Probability, in its etymological sense,
is by nature doubly subjective. First, because
its essence is rooted in a ``feeling'' of the
``human understanding''\,\cite{Hume_Treatise}.
Second, because its value depends on the information
available at a given moment on a given subject.
Many evaluations are based on the assumed properties
of `things' to behave in some ways rather than in others,
relying on symmetry judgments or on regularities
observed in the past and extended to the future
(at our own risk, hoping not to end up
like the {\em inductivist turkey}).
The question of whether there is
``such a thing as {\em Chance} in the world''\cite{Hume_Enquiry}
(does God play dice?) is not easily settled,
but whatever the answer is,
 ``our ignorance of the real cause of any event has the same influence 
on [our] understanding.''\,\cite{Hume_Enquiry} So, at least for
pragmatic convenience, we can assign to `things' {\em propensities},
seen as  parameters of our models of
reality, just like physics quantities.
And they might change with time, as other parameters do.
Furthermore, it is a matter of fact that, besides text book stereotyped cases, propensities
are usually uncertain and we have to learn about them by doing
experiments and framing the observations in a (probabilistic)
{\em causal model}.
The key tool to perform the so-called probabilistic 
inversion is Bayes rule and such models of reality go under the name of
{\em Bayesian networks}, in which probabilities
are attached to all uncertain quantities (possible observations, parameters
and hyper-parameters, which might have different meanings,
including that of propensity and of degree of belief, as when we
model the degree of reliability of a witness in Forensic Science
applications). 
Predictions are then  made by averaging values
of propensities
with weights equal to the probabilities we assign to each of them.

In this paper I have outlined this (in my opinion)
natural way of reasoning, which was that of the founding fathers
of probability theory, with a toy experiment. 
Then, once we have mustered up the
courage to talk about probabilities of probabilities, as shyly done
nowadays by many, we extend them to related concepts, like odds
and Bayes factors.

I would like to end reminding de Finetti's
{\sl ``Probability does not exist''} (in the things),
adding that {\em ``propensity might, but it is in most cases
uncertain and it can change with time}.''

{\small
\begin{flushright}
{\sl \small ``To make progress in understanding,}\\
{\sl \small   we must remain modest and allow that we do not know. }\\
{\sl \small   Nothing is certain or proved beyond all doubt.}\\
{\sl  \ldots} \\
{\sl \small   The  statements of science are not of what is true and what is not true,}\\
{\sl \small  but statements of what is known to different 
             degrees of certainty.''}\\
(Richard Feynman)\\
\end{flushright}
}

\newpage
\section*{Acknowledgements}
This work was partially supported by a grant from Simons Foundation
which allowed me a stimulating working environment during
 my visit at the Isaac Newton Institute
of Cambridge, UK (EPSRC grant EP/K032208/1).
% at the Isaac Newton Institute
%of Cambridge, UK (EPSRC grant EP/K032208/1).
It is a pleasure to thank Dino Esposito for the
long discussions on several of the important points touched upon here and that
we first wrote in our dialogue book~\cite{cep}, and for his patient
and accurate review of the manuscript, which has also benefitted of
comments by Norman Fenton and Geert Verdoolaege.

\end{document}